# Efficient estimators : the use of neural networks to construct pseudo panels

Marie Cottrell[*] and Patrice Gaubert[**]

* MATISSE-SAMOS, Université de Paris I Panthéon-Sorbonne
90, rue de Tolbiac, 75634 Paris 13, France
cottrell@univ-paris1.fr

** MATISSE-SAMOS and LEMMA, Université du Littoral
gaubert@univ-littoral.fr

*Abstract* --- Pseudo panels constituted with repeated cross-sections are good substitutes to true panel data. But individuals grouped in a cohort are not the same for successive periods, and it results in a measurement error and inconsistent estimators.
The solution is to constitute cohorts of large numbers of individuals but as homogeneous as possible.
This paper explains a new way to do this: by using a self-organizing map, whose properties are well suited to achieve these objectives.
It is applied to a set of Canadian surveys, in order to estimate income elasticities for 18 consumption functions..

## 1   Introduction

The need for panel data is synthesized by Baltagi [1] who enumerates the main advantages as follow:
- to control for the individual heterogeneity
- to obtain more information, more variability, less collinearity between variables, more degrees of freedom and more efficiency
- to study more precisely the dynamics of adjustment
- to identify and to measure some effects which are not detectable when using cross-sections data or time series alone

and few other more. According to Verbeek [18] estimators based on panel data are more accurate and more robust for an incomplete specification.
The lack of this kind of data or their inadequacy to be used for specific studies (for instance the PSID and consumption behaviour) leads to the construction of pseudo panels from repeated cross sections. It has been showed by Deaton [7] that the estimators obtained this way possess the same properties as those obtained from true panel data.
Nevertheless, some specific problems arise with construction of pseudo panels, as a result of the grouping of individuals to constitute the cohorts. The fact that the individuals are not the same in two successive observations of the same cohort result in inconsistent estimators. It may be analyzed as a problem of measurement error. Deaton has proposed, in his seminal paper, a treatment of the resulting bias. Verbeek and Nijman [19], analyzing carefully the different aspects of this problem, specially the asymptotic properties, show that this solution is practically inappropriate and define the conditions leading to consistent estimators.

The paper is organized as follows.
In a first part the factors influencing the properties of the estimators are presented according to the results obtained by Verbeek et al. to define the conditions to be respected in the construction of the cohorts of a pseudo panel.
In a second part a Kohonen map is used on a set of Canadian surveys and the cohorts obtained are studied in their main characteristics. The first one is the typology of consumption behaviour obtained, which explains why these cohorts are well suited according to the conditions defined by Verbeek et al.
The last part shows how to use these cohorts in the estimation of a demand system according to the AIDS specification. The income elasticities obtained for 18 consumption functions are computed. Two of them are compared with those computed with a true panel (PSID) on a very similar period of time.

# 2 The need for panel data and pseudo panels : properties and problems

## 2.1 Panel data estimators

(1) Definitions. The most important aim using panel data is the correction of the endogeneity bias linked to the heterogeneity of individual behaviors. In the case of consumption behavior it has been shown (Gaubert [13], for instance) that very different behaviors may be identified in a given population, induced by factors most of the time not present in the estimated equation, resulting in biased estimators.

These properties are presented with the most simple model

$$y_{it} = a + bx_{it} + u_{it}$$
$$i = 1,...,N; t = 1,...,T$$

for an observation $i$ at the period $t$. The term $u_{it}$ may be specified using two ways:

- Model I (fixed individual effect)
$$u_{it} = m_i + n_{it}$$
where $\mu_i$ is a non-random non-observable characteristic and $v_{it}$ a random variable with the usual properties ($n_{it} : iid(0, s_n^2)$).

- Model II (random individual effect):
$$u_{it} = m_i + n_{it}$$
with the hypotheses
$$m_i : iid(0, s_m^2), n_{it} : iid(0, s_n^2)$$
and the $\mu_i$ and the $v_{it}$ are independent.

Two tests are used, one to verify the existence of individual effects (Fisher type), and one to choose the correct model between I and II (Hausman type, see Hausman [9]).

(2) Pseudo panels. Due to the lack of true panel data, Deaton [7] has demonstrated that it is possible to use repeated cross-sections of a population (with completely different individuals from one sample to the next one) to construct pseudo panels and obtain estimators with the same properties as those obtained with true panel data.

(3) According to Deaton cohorts are constituted by grouping the individuals of a survey using a variable which does not change across time, like the birth date, to ensure that this characteristic may be used to group the same subpopulation on repeated cross-sections. This variable may be combined with a few other ones in order to obtain a better definition of the cohorts. Each cohort is created using the same characteristics with a set of cross-sections. This leads to the time dimension of this new unit obtained by averaging the different variables over all the individuals belonging to the cohort.

## 2.2 Problems

Of course the individuals are not the same in a given cohort $c$ created in two consecutive surveys. Even the number of individuals is changing, so the "fixed" effect obtained by averaging $\mu_i$ over the individuals is changing with time for the same cohort $c$.

Moreover, this effect is correlated with the $x_{it}$ in a major part of the economic relations represented with this kind of model, resulting in the inconsistency of random effect estimators.

Two problems have to be distinguished (see Deaton [7] and essentially Verbeek [19]).

(1) A measurement error causes the inconsistency of the estimators. Deaton already explained it but defining a clearly inappropriate corrected procedure of estimation: the required asymptotic property of the estimators in the time dimension is, obviously, never encountered.

(2) A loss of efficiency is due to the grouping of individuals: the estimation on grouped data leads to a loss of efficiency, as it is known since Cramer [6], and the question of how to constitute efficiently the groups has been extensively treated by Haitovsky [14] for a general purpose and by various authors later in specific cases.

According to these authors the reduction of both problems may be obtain being very careful when constituting the cohorts.

- The two objectives pursued in order to obtain consistent and efficient estimators are important but the construction of pseudo panels is one of the cases where there is a trade-off between them.
- Minimizing the within variance of cohort means, relatively to the variance of measurement error, is achieved by grouping very similar individuals in each cohort, and reducing the number of individuals accepted in each one. The definition of numerous cohorts is in favour of more precise estimators, but the small numbers of individuals in each one implies that the computed mean is a poor estimator of the true mean of the population.
- Conversely, the reduction of measurement error obtained with large cohorts results in heterogeneous ones, and a few number of units on which the model have to be estimated.
- Then the solution is to use a set of variables presenting a good adequacy with the studied phenomenon in order to constitute homogeneous cohorts relatively to the significant variables used to describe it.
- While doing this, and depending on the total number of individuals surveyed in the repeated cross-sections, it is necessary to control for the number of cohorts and the number of individuals grouped in each one.
- The constitution of cohorts with the simple cross-classification, obtained using a small number of qualitative variables like age (date of birth), education and so on, is certainly hazardous relatively to the precise sources of bias.

A specific technique seems to be more appropriate.

## 3 Cohorts defined with a Kohonen map

The aim of this application is to construct cohorts having the following properties:
- to be defined using factors quite stable over time in order to link reasonably the successive observations of each cohort
- to be strongly homogeneous relatively to the phenomenon studied (here the consumption behaviour) and as different as possible between them to obtain precise estimators
- to include a number of individuals large enough to allow the use of asymptotic reasoning on the obtained estimators.

### 3.1 The data

We use 3 Canadian surveys performed in 1982, 1986 and 1992 on, respectively, 10936, 9915 and 9475 households.

Consumption expenditures are available for 18 functions together with many socio-economic variables about the household (total income, size, region of residence, tenure status) or the head of family (age, level of education, occupational status, etc.).

A filter is used to exclude a few number of outliers, such as households with negative income and more generally, people with negative consumptions.

The structure of the household consumption (budget shares) and some variables used in the model are presented in Table 1 for the 1986 survey

Table 1 - Descriptive statistics: consumption functions and socioeconomic variables (1986 Survey)

| Functions (budget shares) | Mean | Std dev. |
|---|---|---|
| *Alcohol/Tobacco* | 0.041 | 0.048 |
| *Food at home* | 0.151 | 0.080 |
| *Food away from home* | 0.040 | 0.045 |
| *Housing maintenance* | 0.048 | 0.044 |
| *Communications* | 0.020 | 0.016 |
| *Others (financial costs)* | 0.055 | 0.041 |
| *Transfers* | 0.041 | 0.066 |
| *Education* | 0.014 | 0.030 |
| *Clothing* | 0.070 | 0.043 |
| *Housing* | 0.178 | 0.113 |
| *Leisure* | 0.063 | 0.057 |
| *Furniture* | 0.041 | 0.043 |
| *Health* | 0.023 | 0.025 |
| *Security* | 0.048 | 0.046 |
| *Personal care* | 0.025 | 0.014 |
| *Personal transport* | 0.074 | 0.055 |
| *Public transport* | 0.015 | 0.023 |
| *Vehicles* | 0.049 | 0.108 |
| Total expenditures<br>Age<br>Size of the household (Oxford) | 28291.789<br>47.724<br>2.135 | 16894.090<br>16.258<br>0.933 |

Number of observations: 9606

### 3.2 Variables used to make the classification

(1) To construct cohorts as homogeneously as possible, relatively to the model we have to estimate, and at the same time, to obtain cohorts as different as possible between them, the best indication is to use the significant variables of these behaviours: a cohort constituted mainly with households sharing more or less the same order of preferences seems to present the first property required to construct correctly the pseudo panel. So, the principal variables used are the budget shares defining the consumption behaviour of each household.
- This structure is not as independent over time as is the date of birth (in Deaton's cohorts): the shares are varying over the period of observation, according to the Engel laws, but
- the period of time is not too long to make the hypothesis that these behaviours are only slightly varying, as it is verified when the different resulting cohorts are analysed
- according to the studies cited above about the efficient method to group data, the explanatory variables of the model are included in the input data space of the algorithm[1], with a special treatment of age.
- This means that the algorithm is used to construct groups using quantitative variables, the qualitative ones being used only to interpret the groups obtained.

(2) Age of the individuals has to be treated in a specific way: if it is used to construct the classification like the other variables, the result is that the classes obtained are very homogeneous considering all the variables used, including age. This is a main difference if we compare these cohorts with those produced with the Deaton's method. The problem is that this group is constituted of individuals who share the same consumption behaviour at the same age. So a part of the dynamic process of consumption, the fact that people reach a level of consumption in some good at different steps of their life cycle, is concealed. Constraining the group to have a common age does presumably produce a kind of cohort with an averaged behaviour leading to measures similar to the ones obtained on a simple cross-section.

- To avoid this, dummy variables are created to represent classes of age and the limits of these classes are corrected[2] in the successive cross-sections to conform with the idea of cohorts as defined by Deaton.
- A more simple treatment of age is tried simultaneously: controlling for age varying only with the inclusion of this variable in the specification of the estimated model.

### 3.3 Classification and constitution of the cohorts

The map constructed is a grid of 64 nodes.
The performance of the Kohonen algorithm to reveal very differentiated behaviours of consumption has been presented in a former study (Cottrell et al., 2000).

(1) Some examples are presented to verify the quality of the classification to create groups which have very different behaviours. These behaviours may be briefly described, and related to qualitative variables characterizing the household and its composition:
- for the whole sample the main functions are "food at home" (15 %) and "housing" (18 %), then "personal transport" (8 %), "clothing" (6.6 %) and "leisure activities" (6.2 %); the shares of the 13 other functions are between 2 and 4 %.
- class 1 has a consumption dominated by vehicles expenditures, it is constituted with households of two adults older and receiving higher resources than the population average
- class 4 is made of households of two adults older than the preceding ones and with a high level of expenditures devoted to the transfers
- class 13, with households of one adult of middle age and one child, is characterized by a high level of "Food away from home" and "Housing"
- class 22 with 2 adults older than the average is dominated by "Health" expenditures
- classes 24 and 32 with one or two old and poor adults are devoting most of their incomes to "Food at home" and "Housing"
- conversely, class 57, with households of two adults and two children, younger and richer than the average, have a consumption behaviour dominated by "Education" expenditures.

These are only a few examples to enlighten on the quality of the output produced by the neural algorithm.

(2) More precisely the differences between consumption behaviours may be measured in

---
[1] Besides of that the model estimated include a time variable in order to capture the effect of changing environment as well as changing tastes.
[2] 4 years more in 1986 and 6 more in 1992.

order to verify that the main objective defined when constructing cohorts is reached.
- Distances between classes. The Mahalanobis distances between the nodes may be computed using the code vectors at the end of the iterative process.

They may be represented on a grid similar to the Kohonen map[3] using polygons which express the distance between a node and its 8 closest neighbours: the more the polygon is far from the contours of the cell the more distant is the node from the corresponding neighbours.

This map shows clearly the significant differences between the nodes, even with the use of neighbouring during the whole process of construction of the map[4].
- Within and total variances. Another way to measure the quality of this classification is to measure the share of total variance computed for the whole sample resulting in within variance when the cohorts are defined: the smaller this variance relatively to total variance, the more homogeneous the cohorts.

A non-parametric test may be computed, the Wilks test, on the distribution of the budget shares over the 64 groups obtained.

A comparison between this measure obtained with the Kohonen cohorts, and the one produce by a Deaton-like construction (age combined with level of education and region of residence) on the same set of Canadian surveys and with the same number of groups is presented (Table 2)

Table 2 - Share of within variance relatively to total variance with neural and Deaton-like cohorts

|  | SOM | Deaton-like |
|---|---|---|
| Alcohol-Tobacco | 48.95 | 94.59 |
| Food at home | 51.12 | 85.55 |
| Food away from home | 50.03 | 96.49 |
| Housing maintenance | 47.20 | 86.59 |
| Communication | 77.22 | 94.68 |
| Others | 47.10 | 98.13 |
| Transfers | 34.25 | 87.19 |
| Education | 25.26 | 89.81 |
| Clothing | 55.46 | 94.66 |
| Housing | 45.35 | 91.87 |
| Leisure | 43.99 | 96.04 |
| Furniture | 48.17 | 98.52 |
| Health | 68.11 | 96.49 |
| Security | 46.82 | 77.72 |
| Personal care | 70.65 | 97.84 |
| Personal transport | 60.93 | 95.29 |
| Public transport | 49.27 | 95.77 |
| Vehicles | 16.80 | 98.71 |
| Wilks Lambda | 0.00000623 | 0.3584 |
| F | 145.27 | 9.49 |

The neural classification dramatically reduces the within variance, compared with the classical construction. It appears that the latter has quite no relationship with the phenomenon of interest: only 5 functions on 18 show a within variance lower than 90%.

Conversely, for the neural classification, with the exception of 4 functions known to put together heterogeneous goods or services (communication, health, personal care and personal transport), the within variance has a share lower than 50%.

(3) We have to check now the other constraints imposed to the cohorts in order to reduce most of the measurement error and obtain efficient estimators.

It is possible to represent a portion of the network constructed, summarizing the contents of the classes with the number of observations which belong to each survey.

If the classes are numbered from the top left corner (1) to the bottom right corner (64) going from top to bottom and from left to right, the following extract shows the number of individuals of each survey gathered in one of these four classes:

| C4 | C12 |
|---|---|
| 1982 survey: 181 | 1982 survey: 109 |
| 1986 survey: 158 | 1986 survey: 123 |
| 1992 survey: 139 | 1992 survey: 172 |
| C5 | C13 |
| 1982 survey: 136 | 1982 survey: 165 |
| 1986 survey: 121 | 1986 survey:  89 |
| 1992 survey: 118 | 1992 survey: 114 |

For instance class 4 is constituted of two adults older than the average, they have a high budget share for "Transfers". This constitutes a cohort, $C_i$ thereafter and $C_{it}$ for the observation of this cohort at the period $t$, of 181 individuals observed in 1982, 158 individuals in 1986 and 139 in 1992. These 3 groups have something common: a specific consumption behaviour which produces their assignment to this class which is closer to their own behaviour than any other.

---
[3] See Appendix Fig. 2.
[4] Except the last iteration which is usually executed with a neighbouring distance reduced to one : only the winner node have its code vector adapted.

The size of each class varies across the map, but there are only 12 classes with less than 300 individuals, meaning that these cohorts gather less than 100 individuals for one of the surveys or more. Only 7 cohorts have between 150 and 225 observations. All the others are more numerous, with a number of observations in each $C_{it}$ greater than 100.

According to the computations produced by Verbeek et al. about the size of the bias, these cohorts seem to present the right properties to obtain consistent and efficient estimators.

In the following we work with this pseudo panel which consists of 64 "statistical" individuals measured three times.

### 3.4 Application : consumption functions

(1) Our pseudo panel of consumption expenditures produced by the neural treatment of 3 cross-sections leads to the estimation of demand functions using an AIDS specification, as it has been defined by Deaton and Muellbauer [8].

According to Banks et al. [2] quadratic terms are added in order to capture some non-linearities which appear to be significant for some of the $j$ functions (QUAIDS).

Differentiated prices for each survey are not available, and the price variables have to be removed from the model. The effect of changing prices as well as the effect of changing environment will be taken into account through the use of a time effect.

Due to the classical measurement error attached to the household's incomes, the total expenditure has been substituted as an instrument. It is well known that this instrument may be itself affected by a measurement error, being the sum of items diversely concerned by this type of error. The hypothesis used here is that this error is cancelled by grouping the data into cohorts.

As is usual for this type of estimation, 2 control variables are added: the age of the household's head and the family equivalized size.

As a result, the equation[5] to be estimated (QUAIDS) is

$$w_{it} = a_0 + b_1 \log y_{it} + b_2 (\log y_{it})^2 + b_3 \log age_{it}$$
$$+ b_4 (\log age_{it})^2 + b_5 \log size_{it} + b_6 (\log size_{it})^2$$
$$+ d \, year_t + \mu_i + n_{it}$$

where $w$ is the budget share (of one function) for the individual $i$ (the mean of the cohort $C_i$) at the period $t$, $y$ is the total expenditure, *age* the age of the household's head, *size* the size of the household using the Oxford scale; *year* is time dummy and $\mu$ is the individual effect.

For the AIDS specification, the $\beta_2$ term is removed.

(2) Constructing the variables of the statistical individuals[6]. We estimate this equation at the cohorts' level not on the individual values: for each class produced by the classification we compute the mean of every variable used in the model.

Because the dependent variable is the share of expenditure in one good relatively to total expenditure, the computation of the mean of each other variable has to be weighted using a factor

$$g_{ht} = \frac{y_{ht}}{\sum_{h \in C_{it}} y_{ht}}, \text{ so that } w_{it} = \sum_{h \in C_{it}} g_{ht} w_{ht}$$

All the variables used in the model are computed this way, as a weighted mean of the values measured at the individual level in each cohort.

The heteroscedasticity introduced by this construction has to be removed: this is done by pre-multiplying each variable by the inverse of the square root of the factor inflating the residual variance, that is $\frac{1}{\sqrt{\left(\sum_{h \in C_{it}} (g_{ht})^2\right)}}$.

The model is estimated for the observations constituted by these transformed values computed for the 64 cohorts observed over 3 periods.

---

[5] For each function, the subscript of which being omitted in order to facilitate the reading.

[6] See a detailed presentation of the transformation in Cardoso et al. [3]

The AIDS and QUAIDS specifications are systematically used and tested to identify the one in adequacy with the data. At the same time, models I and II (fixed or random individual effects) are successively tested, even if the fixed effect specification is generally preferred due to the correlation between the effect and the explanatory variables. The Hausman test is used to indicate the more convenient.

(3) The elasticities. The elasticities computed from the estimated parameters (Table 3) present a very good level of accuracy

Table 3 - Total expenditures elasticities

| Functions | Class. with age as dummies | | Age only in the model | |
|---|---|---|---|---|
| | Elasticity | Student | Elasticity | Student |
| Alcohol-Tobacco | 0.730 | 7.686 | 0.879 | 11.459 |
| Food at home | 0.474 | 10.170 | 0.419 | 9.153 |
| Food away from home | 1.275 | 17.118 | 1.185 | 13.709 |
| Housing maint. | 0.636 | 8.006 | 0.384 | 3.992 |
| Communication | 0.846 | 13.171 | 0.779 | 11.701 |
| Others | 0.974 | 14.234 | 1.075 | 16.568 |
| Transfers | 1.316 | 14.475 | 1.410 | 10.671 |
| Education | 1.258 | 19.353 | 1.460 | 14.104 |
| Clothing | 0.963 | 17.182 | 1.004 | 13.214 |
| Housing | 0.905 | 15.674 | 0.990 | 13.183 |
| Leisure | 1.258 | 30.941 | 1.246 | 19.584 |
| Furniture | 1.005 | 14.610 | 0.948 | 10.040 |
| Health | 1.045 | 13.754 | 1.142 | 11.281 |
| Security | 1.322 | 27.175 | 1.262 | 11.645 |
| Personal care | 0.848 | 13.406 | 0.855 | 16.396 |
| Personal transport | 0.908 | 9.554 | 0.857 | 7.806 |
| Public transport | 1.091 | 10.483 | 0.768 | 5.949 |
| Vehicles | 1.898 | 8.988 | 1.919 | 7.661 |

For 10 functions over 18, the adequate specification is AIDS, and for 6 functions the Hausman test rejects the error component form.

These results are consistent with the general ideas on necessary goods (elasticity significantly lower than 1) and dynamic goods (goods that budget share grows with the income).

For two functions, "food at home" and "food away from home", a comparison may be tempted: income elasticities have been obtained using the same specification for a similar period of time (1985-87) using true panel data on American families (PSID)[7]. The numbers are very close: respectively 0.24 and 0.80 with PSID data, 0.23 and 0.89 with this pseudo panel.

## 4    Conclusions

The construction of a pseudo panel from repeated cross-sections using a neural network like the Self-Organizing Map appears to be a means to overcome the major drawbacks attached to the classical pseudo panels.

It produces cohorts with a great homogeneity relatively to the phenomenon studied, depending on the variables chosen to constitute the input data. As the principle of this technique is to transform multidimensional data into a structure compressing this information while preserving the essential, that is the initial topology, the result corresponds to the aim if the variables used are pertinent.

The only limit to the quality of the pseudo panel obtained is the number of individuals available in each survey: in order to use asymptotic reasoning to evaluate the estimators and to obtain accurate estimations the number of cohorts must be greater than 50 and the size of each cohort has to be at least 100 observations. This determines the minimum size of the surveys, considering that the algorithm produces reasonably balanced classifications, but there is some variation in size between the classes produced.

The combination of a set of pertinent variables in the input data of the algorithm gives the opportunity to use qualitative variables to construct the cohorts, even if the estimated model is not suited to include them with a fixed effect specification.

**Acknowledgements**

The data is obtained from Statistics Canada through Pr. Simon Langlois, Université Laval (Québec).
The authors thank François Gardes for numerous suggestions and comments, and participants of Paris I seminar, LEMMA seminar and 9th International Conference of ACSEG (2002) for their comments.

---

[7] See Gardes et al. [10] : total expenditures elasticities have to be multiplied by a factor of 0.7, the income elasticity of total expenditures.

# APPENDIX
# THE CLASSIFICATION [8]

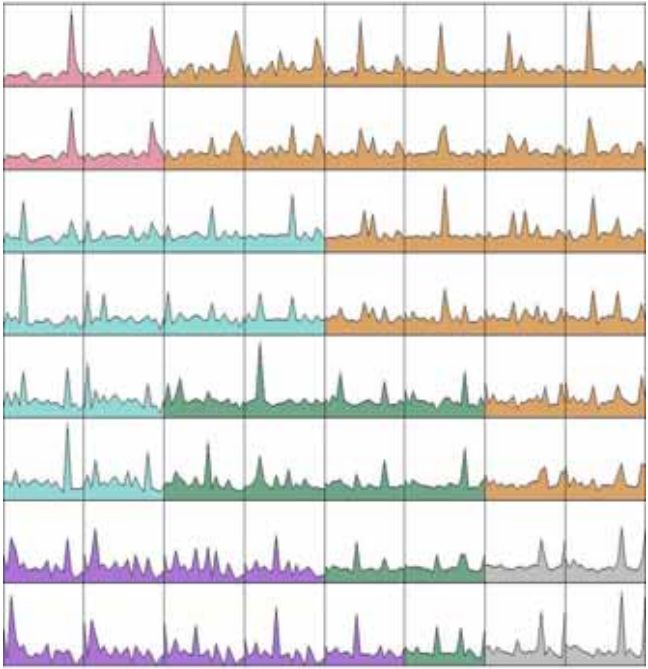

Fig. 1. The Kohonen Map: a representation of the code vectors.

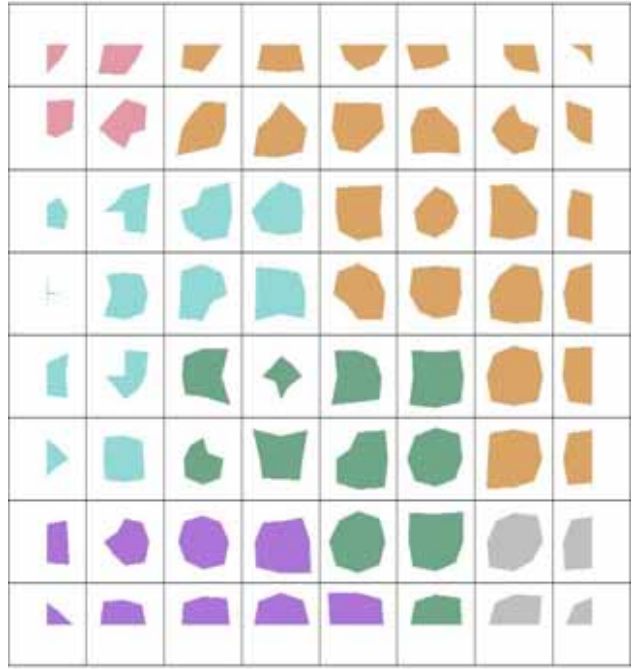

Fig. 2. The Kohonen Map: the distances between the classes.

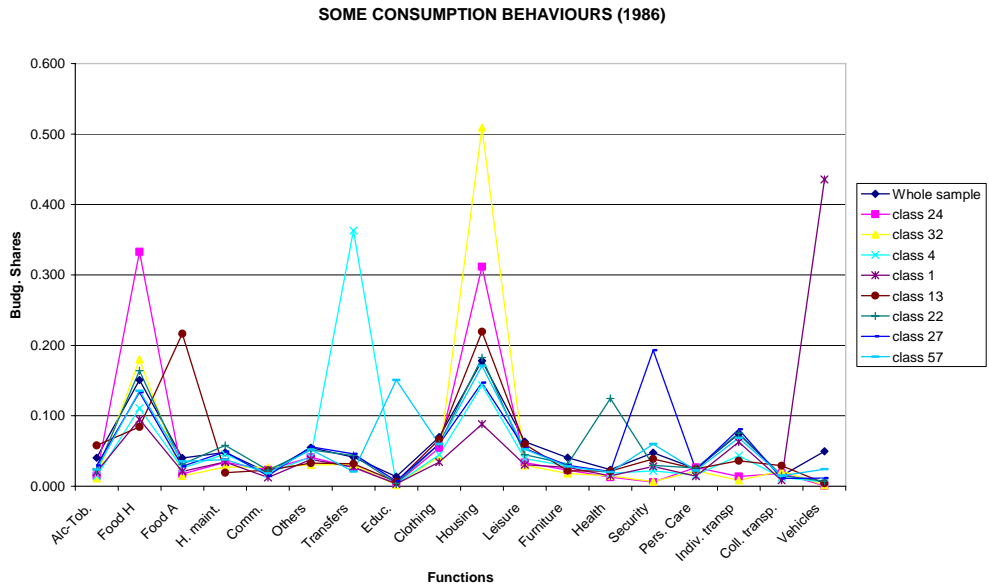

Fig. 3. Budget shares of the whole sample and a selection of the 64 cohorts obtained

---

[8] The programs used to construct the network and to realize some further statistical treatments may be obtained at http://samos.univ-paris1.fr .